\documentstyle{article}

\input amssym.def
\input amssym.tex

\newtheorem{theorem}{{Theorem}}[section]
\newtheorem{hermann}[theorem]{{Frobenius' Theorem for finitely
generated distributions}}
\newtheorem{integrability theorem}[theorem]{{Integrability Theorem}}
\newtheorem{rough structure theorem}[theorem]{{Rough Structure Theorem}}
\newtheorem{nother}[theorem]{{Noetherian Theorem}}
\newtheorem{lemma}[theorem]{{Lemma}}
\newtheorem{corollary}[theorem]{{Corollary}}
\newtheorem{remark}[theorem]{{Remark}}
\newenvironment{proof}{{\it Proof.}}{\hfill$\diamondsuit$\medskip}

\begin{document}

\title{On Gromov's theory of rigid transformation groups: A dual approach}
\author{Abdelghani Zeghib}
\date{}
\maketitle

\begin{abstract}

Geometric problems are usually formulated by means of   (exterior) differential
systems. In this theory, one enriches the system by adding
algebraic and differential constraints, and then looks for regular solutions.
 Here we adopt a dual approach, which consists to enrich
a plane field, as this is often practised in control theory,
 by adding brackets of the vector fields tangent to it, and then, look for
singular solutions of the  obtained distribution.  We apply
this to the isometry problem of rigid geometric structures.


\end{abstract}



\medskip
\paragraph{Content.}In \S \ref{Control}, we exhibit a natural  class of
plane fields,
for
which the accessibility behaviour,  as studied in control theory,
possesses, essentially, the same nice properties, as in the analytic case.

In \S \ref{Geometry}, we observe that there is a control
theory approach to  the local isometry problem of affine manifolds
(e.g. pseudo-Riemannian manifolds), which is dualy equivalent, to
the usual differential systems (i.e. partial
differential relations) approach. We then apply the results
of \S \ref{Control}, to deduce a celebrated corollary of
Gromov's theory on rigid transformation groups.

In fact, the developments of \S \ref{Control}, suggest how to proceed  in
order to
recover  essentially  ``most of the   theory'', together with  some independent
results. But, we won't follow this, because, our primary
goal here is to be as elementary as possible.


\section{Control theory}
\label{Control}

Let $P$ be a smooth plane field of dimension $d$ on a manifold
$N$. From an integrability viewpoint,  there are two extremal cases,
described by the classical  Frobenius and Chow's theorems,
which concern the completely integrable and
absolutely
 non-integrable cases, respectively.

 Let
{\Large{$\chi$}}$(N)$ be the Lie algebra of $C^\infty$ vector fields
of $N$, and denote by
$\cal G$  the Lie subalgebra generated by
smooth vector fields everywhere
tangent to
$P$.
Let $G$ be   the ``evaluation plane field'',  $G(x) =
\{ X(x), X \in {\cal G} \}$ (this is
not necessarily a continuous plane field).

Frobenius' Theorem  states that in  the ``degenerate case'' where $P$ is
involutive, that is
$G= P$, then through each point of $N$,
passes a {\it leaf} of $P$, that is a submanifold of dimension
$d$ (the same as that of  $P$) which is (everywhere) tangent to $P$. In
contrast,
in the ``generic case'',
when $G = TM$, Chow's Theorem says that any pair of points can
be joined by a curve tangent to $P$. However (unfortunately),
it is the intermediate
(non-generic and non-degenerate) situation
that one usually meets in geometric and differential problems.


\paragraph{Integrability and infinitesimal integrability domains.}
In searching leaves, let's  ``naively''  introduce the { \it integrability
domain}
${\cal D}$ as the set of points of $N$, through which passes a (germ) of a
leaf of $P$.

This set may behave very badly, for instance, it is
not a priori closed.
For this, let's introduce its infinitesimal variant,
the ``involutivity domain'',
${\cal D}^\infty= \{x \in N \; / G(x) = P(x) \}$

We call ${\cal D}^\infty$ the {\it infinitesimal integrability domain}
of $P$.
Clearly, ${\cal D}^\infty$ is closed and contains
${\cal D}$.

Along ${\cal D}^\infty$, the Frobenius
condition is satisfied, and so, one may hope to
find leaves through each of its points, that is
${\cal D} = {\cal D}^\infty$. However,  ${\cal D}^\infty$ is
not,  a priori,  a manifold, and we do not yet know a fractal
Frobenius' Theorem. Worse, even if we assume ${\cal D}^\infty$
is a submanifold, it is not clear that $P$ is tangent to it!

In the analytic case,  everything works well, and there
are many ways leading to the equality
${\cal D} = {\cal D}^\infty$.

For instance ${\cal D}^\infty$ is an analytic set,
and may be thought out as being a submanifold, and so
in order to apply Frobenius' Theorem to the restriction
$P \vert {\cal D}^\infty$, one
 just has  to show that $P$ is tangent to
${\cal D}^\infty$.

\paragraph{Distributions.} However, the  most consistent
 approach to this problem
is a generalization of Frobenius' Theorem in another
direction, that of (singular) {\it distributions }
(the singularity is topological and not differential).
Recall that a
$C^\infty$ distribution $\Delta$ on $N$
 is a $C^\infty(N)$-submodule of {\Large{$\chi$}}$(N)$,
the space of $C^\infty $
vector fields on $N$.

For example, to a smooth plane field $P$
is associated the distribution of vector fields tangent to it.
Conversely, to a distribution $\Delta$, one defines its
 ``evaluation
plane field''  by $\Delta(x) = \{X(x) / \; X \in \Delta\}$. In general, this
 determines a discontinuous plane field (i.e. a plane field with
 non-constant dimension). One calls a distribution {\it regular}
if its ``evaluation plane field''  has constant dimension.

A distribution is called {\it involutive} if it is a
Lie  subalgebra of {\Large{$\chi$}}$(N)$. Any distribution generates
an involutive distribution, this is the advantage
of generalizing plane fields to distributions, since
the involutive distributions generated by plane fields,
are not plane fields in general, i.e. they aren't necessarily
 regular.

\paragraph{The integrability problem.}
A {\it leaf}  of a distribution
$\Delta$
is a submanifold $S$  such that along $S$ the tangent space of
$S$,  {\it coincides} with the evaluation of
$\Delta$. The distribution is called
 {\it integrable} if leaves
exist everywhere, i.e. any $x \in N$ belongs to a leaf.


In particular,
if
the involutive distribution ${\cal G}$ generated by
a plane field
$P$  is integrable, then we have in particular  the equality
${\cal D} = {\cal D}^\infty$. Indeed, if $x \in
{\cal D}^\infty$, then its ${\cal G}$-leaf, is
a leaf of $P$. (In the generaic, but non-interesting case, ${\cal G}$
is integrable, and ${\cal D}^\infty = \emptyset$).

Obviously, an integrable distribution is involutive.
The converse is not true in general (see \ref{counter.example} for
counter-examples).
The integrability
problem consists of finding
conditions so that involutive implies integrable.
For instance, Frobenius' Theorem says nothing but that   regular
involutive
 distributions  are integrable.

\paragraph{Finitely generated distributions.}  A distribution $\Delta$
is
called  {\it locally finitely generated}, if for any $x \in N$, there is
a neighborhood $U$ of $x$, and a finite family
$V_1, \ldots, V_l$ of
vector fields of $\Delta$, such that,   on
$U$, any $V \in \Delta$ can be written as
$V = \Sigma g_iV_i$, where $g_i \in C^\infty(U)$.
Regular distributions  are locally finitely {\it freely} generated.
Conversely, it
seems that a suitable ``blowing up'' manipulation transforms
locally finitely generated distribution, to
regular ones. Anyway, Frobenius' Theorem is valid in this context.

\begin{hermann}  \label{hermann}
(R. Hermann \cite{Her}, see \S  \protect{\ref{proof.hermann}} for an
outline of proof).
A locally finitely generated involutive distribution is integrable.
In particular, let $P$ be a smooth plane field such
that its associated involutive distribution is finitely generated. Then
${\cal D} = {\cal D}^\infty$.
\end{hermann}

\paragraph{Partially algebraic vector fields.}
The theorem above applies in the analytic case, thanks to standard
Noetherian facts. We are now going to
extend the applicability of the theorem above to
a partially analytic, in fact partially algebraic,  situation.
The starting point is to consider
{\it partially algebraic} vector fields
on ${\bf R}^n \times {\bf R}^m$. They are $C^\infty$
vector fields of the form: $(x, u) \in
{\bf R}^n \times {\bf R}^m \to
(R(x, u), Q(x, u))$ such that, for $x$
 fixed,  $R(x, u)$ and $Q(x, u)$ are polynomials.

In other words, partially algebraic  vector fields
are mapping: $ {\bf R}^n \times {\bf R}^m \to
{\bf R}^n \times {\bf R}^m $, with
co-ordinates in the ring
$C^\infty( {\bf R}^n) [X_1, \ldots,  X_m]$.

Observe that the bracket of two partially
algebraic vector fields is a partially
algebraic vector field. That is,  partially algebraic
vector fields form a Lie subalgebra.

 Let $\Phi$
 be a partially
linear (local) diffeomorphism of
${\bf R}^n \times {\bf R}^m $,
 that is $\Phi $
has the form $\Phi(x, u) = (f(x), A_x(u) )$, where
 $f: U \to U^\prime$ is a local diffeomorphism
of ${\bf R}^n$, and $A: x \in U \to A_x \in GL(m)$
is a $C^\infty$ mapping.

Observe that partially linear diffeomorphisms
preserve the space of partially algebraic
vector fields. (Here one can
also consider partially
polynomial diffeomorphisms, but for the sake of
simplicity, we restrict ourselves to
the partially linear case).

\paragraph{Fiberwise algebraic vector fields on vector bundles.}
 Suppose that $N \to B$ is a vectorbundle.
What precedes, allows us to define  {\it fiberwise
algebraic vector fields} on $N$.
They form a Lie
subalgebra.

One can also define {\it fiberwise algebraic
plane fields}  and {\it fiberwise algebraic distributions}.

The involutive distribution
generated by a fiberwise algebraic distribution is
fiberwise algebraic.

One can also define fiberwise algebraic functions, and then
{\it fiberwise algebraic sets}, as zero loci of systems
of fiberwise algebraic functions.

\begin{integrability theorem} \label{integrability}
Let $P$ be a fiberwise
algebraic plane field on a vector bundle $\pi: N \to B$.
Then, there is an open dense set $U \subset B$,
over which ${\cal D}= {\cal D}^\infty$. More
precisely, the involutive distribution generated
by $P$ is integrable on $\pi^{-1}(U)$.

\end{integrability theorem}

\begin{proof} Let ${\cal G}$ be the involutive
distribution generated by $P$.
From the previous discussion, it can be described locally as
an $R$-submodule $I$ of $R^{n+m}$, where
$R= C^\infty({\bf R}^n)[X_1, \ldots, X_m]$.

Following
Theorem \ref{hermann}, it
suffices to show that over an open dense set
$U$ of ${\bf R}^n$, $I$ is locally finitely
generated. This will follow from the
Noetherian Theorem \ref{nother}. The intuitive proof of it,
is that we have a family $\{ I_x ; \; x \in {\bf R}^n\}$ of ${\bf R}[X_1,
\ldots,
X_m]$-submodules of $({\bf R}[X_1, \ldots, X_m])^{n+m}$. Each
$I_x$ is finitely generated and,  in a dense open subset of ${\bf R}^n$,
the cardinality of the generating family of $I_x$  is locally bounded.



\end{proof}

\paragraph{Differential structure of ${\cal  D}$.}  The infinitesimal
integrability  domain
${\cal D}^\infty$ (and hence the integrability domain
${\cal D}$, if we restrict over $U$) is
a  fiberwise algebraic set. Indeed, ${\cal D}^\infty$ is
the set of points where the involutive distribution ${\cal G}$
generated by $P$, has dimension $d$ (that is the dimension of $P$).
Thus, ${\cal D}^\infty =
 \{ x \in N / \; V_1 \wedge \ldots \wedge V_{d+1} = 0 $
 for all $V_1, \ldots V_{d+1}$ elements of ${\cal G} \}$.
 Locally ${\cal D}^\infty $
is the  zero locus of a family of elements
of  $C^\infty ({\bf R}^n) [X_1, \ldots, X_m]$.

The fibers
${\cal D}^\infty_x$ are thus algebraic sets of
${\bf R}^m$.

In fact, fiberwise algebraic sets   enjoy in addition
many basic (that is in dependence on $x \in B$)
regularity properties.

In local co-ordinates, around a point where the
distribution
${\cal G}$ is locally finitely  generated,  ${\cal D}^\infty$ is
the common  zero locus of a finite set
$f_1, \ldots, f_l$ of elements  of $C^\infty({\bf R}^n)
[X_1, \ldots, X_m]$. But, because,
we reason here over ${\bf R}$ (and not ${\bf C}$),
${\cal D}^\infty$ equals the zero locus of a single element
$g = \Sigma f_i^2$.  This
element $g$ may be seen as
 a map $f: {\bf R}^n \to
{\bf R}[X_1, \ldots, X_m]_{\leq k}$, the
space of polynomials of degree
$\leq k$ ($f(x)$ is the restriction
of $g$ to
$\{x\} \times {\bf R}^m$).
More concretely, by definition of
$C^\infty ({\bf R}^n)[X_1, \ldots, X_m]$,
we have a representation
$g(x, X_1, \ldots X_m) = \Sigma_{ \vert I \vert \leq k} g_I (x) X^I$, where
$I$ is a multi-index, then $f(x)$ is
the polynomial with coefficients  $(g_I(x))_{\vert I \vert \leq k}$.

Suppose for example that $f(x)$ has a
unique (real) root $z(x) \in {\bf R}^m$, and thus
${\cal D}^\infty$ is the graph of $z$. Then, $z(x)$ is expressed
``algebraically'' from the coefficients of $f(x)$. Therefore,
${\cal D}^\infty$  is the graph of a very ``tame'' function.

The same idea may be adapted when  $f(x)$ has
infinitely many roots. This may lead to
 a stratified structure of ${\cal D}^\infty$, after removing
singular fibers.
 We will restrict our  investigation  here
to a week regularity aspect, which will follow from
the following general fact.



\begin{lemma}  \label{open.dense}
Let $B$ be a topological space, and  $f:
B \to {\bf R}[X_1, \ldots, X_m]_{\leq k}$, a continuous map which
associates
a polynomial $f(x)$ of degree $\leq k$, to each $x \in B$. Let
$Y = \{x \in B /\; f(x)$ has a (real) root $\}$. Then
$Y$ contains an open dense set of its closure.

\end{lemma}

\begin{proof} Consider the
``universal'' polynomial
 $$ \Phi:  (X_1, \ldots, X_m, p) \in {\bf R}^m \times
{\bf R}[X_1, \ldots, X_m]_{\leq k}
\to  p(X_1, \ldots, X_m) \in {\bf R}$$
($p$ is a polynomial    of
degree $\leq k$
on $(X_1, \ldots, X_m)$).

Consider  the ``universal''
algebraic set $\Phi^{-1}(0)$ determined by $\Phi$. Let $Z $ be
 the projection of
$\Phi^{-1}(0)$ on ${\bf R}[X_1, \ldots, X_m]_{\leq k}$. It isn't a priori
an algebraic set, but,
almost by definition, a {\it semi-algebraic set}.

One fundamental fact about semi-algebraic sets is that they admit
good  stratification (see for example \cite{Ben}). In particular, $Z $
is a finite disjoint union $Z = \cup Z_i$,
where $Z_i$ are {\it locally closed} sets, that is,  there are
open sets  $O_i$
 in
 ${\bf R}[X_1, \ldots, X_m]_{\leq k}$, such that
$Z_i = \overline{Z_i} \cap O_i$.

For the lemma, we may assume that
$Y$ is dense in $B$, we have then to show that
$Y$ contains an open dense set of $B$. By continuity,
$f(B)$ is contained in
$\overline{Z}$ (which also equals
$ \cup \overline{Z_i}$).
 We have, $Y = f^{-1}(Z)$.

If $Z$ itself were locally closed (for example for $m =1$), then $f^{-1}(Z)$
would be open in $B$, and we are done.

We argue as follows in the general case. Let
$F_i = f^{-1}( \overline{Z_i})$, and
$A_i = f^{-1}(Z_i)$.
Then  $A_i$ is open
in $F_i$ (because $Z_i$ is locally closed).

We have,
$B = \cup F_i$.
One firstly observes    that
$\cup int(F_i)$ is dense in $B$, where
$int$ stands for the interior (this is Baire's Theorem,
for {\it finite} union of closed sets, which is true for
all topological  spaces). Next,
since $Y = UA_i$ is dense in $B$, it follows that
 $ U = \cup (A_i \cap int(F_i))$ is dense
in $B$. Moreover $U$ is open in $B$ (since
$A_i$ is open in $F_i$) and is contained in $Y$.

\end{proof}

The discussion before the lemma applies to any
fiberwise algebraic set (like ${\cal D}^\infty$), and therefore leads to
the following result.

\begin{corollary} \label{structure.algebraic} Let ${\cal S}$ be
a fiberwise algebraic set of $N$
 and $Y$ its projection on
$B$. Then $Y$ contains an open dense subset
of its closure $\overline{Y}$.

\end{corollary}

\paragraph{Fiberwise constructible sets.} In view of further  applications,
we need
the following slight generalization of fiberwise algebraic sets. A subset
${\cal S}$ of
$N$ is called {\it fiberwise constructible}  if it can be written as
a difference ${\cal S}_1 - {\cal S}_2$ of two fiberwise algebraic sets
${\cal S}_1$ and ${\cal S}_2$.

Such a set has a  structure as  nice as  that of a fiberwise
algebraic set. Indeed, locally, suppose that ${\cal S}_1$
and ${\cal S}_2$ are respectively defined by $f$ and
$g$ elements of $C^\infty({\bf R}^n) [X_1, \ldots, X_m]$. Then, consider
 the mapping $\phi: (x, X) \in {\bf R}^n \times {\bf R}^m -{\cal S}_2 \to
(x, X, 1/g(x, X)) \in  {\bf R}^n \times {\bf R}^{m+1}$ (here $X = (X_1,
\ldots, X_m) $).
Then, the image $\phi({\cal S}_1 - {\cal S}_2)$ becomes fiberwise
algebraic, since it is defined
by the equations, $X_{m+1} g(x, X) -1 = 0$, and $f(x, X)=0$.

The corollary above is therefore valid for fiberwise constructible sets.

\begin{rough structure theorem} \label{structure}
Let ${\cal S}$ be a fiberwise constructible set of $N$
 and $Y$ its projection on
$B$. Then $Y$ contains an open dense subset
of its closure $\overline{Y}$. In particular if $Y$ is
dense in $B$, then $Y$ contains an open dense
subset of $B$.

\end{rough structure theorem}

\paragraph{Integrability with constraints.} One is sometimes interested
 in  leaves through points in a given subset ${\cal S} \subset N$
(the plane field $P$ is not assumed to be tangent to
${\cal S}$, although this usually happens in  practice).


The following result unifies the two previous theorems \ref{hermann} and
\ref{structure}:

\begin{theorem} \label{constraints}
Let $P$ be a fiberwise algebraic plane field on a vector bundle
$\pi: N \to B$, and ${\cal S}$ a fiberwise
constructible subset of $N$. There is an open dense set $U \subset B$,
over which,
the sets of integrability and infinitesimal
integrability points of $P$ in
${\cal S}$  are equal, that is,
${\cal D} \cap {\cal S} \vert U  =
{\cal D}^\infty \cap {\cal S} \vert U$.

In addition, the projection of ${\cal D} \cap {\cal S} \vert U$ is a closed
(may be empty) subset of  $U$.

\end{theorem}

\begin{proof} Let $U_1$ be an open dense set given by the integrabilty theorem
\ref{hermann}, that is ${\cal D} \vert U_1= {\cal D}^\infty \vert U_1$.
Over $U_1$,
${\cal D} \cap {\cal S}$ is fiberwise constructible. Let $Y_1 \subset U_1$
be its projection, and let $\overline{Y_1}$ be its closure in $U_1$. From
the structure theorem \ref{structure}, there is an open subset $U_2$ of
$U_1$, such that
$Y_1$ contains $\overline{Y_1} \cap   U_2 $, which is in addition dense in
$\overline{Y_1}$. In particular $Y_1 \cap U_2$ is closed in
$U_2$.

We claim that   $U = U_2 \cup (U_1 - \overline{Y_1})$ satisfies the
conditions of
the theorem. Indeed, $U$ is open, and  it is dense in $U_1$ (and hence in
$B$), since
 $\overline{Y_1} \cap   U_2 $ is dense in
 $\overline{Y_1} $. Furthermore, $U \cap Y_1 = U \cap  \overline{Y_1} $ and
hence,
over $U$, the projection of ${\cal D} \cap {\cal S}$ is closed.

\end{proof}



\section{The isometry pseudo-group of an affine connection}
\label{Geometry}

Fiberwise algebraic objects are abundant in geometry.
For instance, a fiberwise algebraic function on  the
cotangent bundle of a  smooth manifold, generates
 a fiberwise algebraic Hamiltonian vector field.
In particular the geodesic flow of a Riemannian metric
is fiberwise algebraic (being seen on
the cotangent as well as on the tangent bundles).

\paragraph{ The tautological geodesic
plane field of an affine manifold.}
More generally, let $(M, \nabla)$ be an affine manifold, that is,  $\nabla$
 is a torsion free connection on $M$ (not necessarily flat).
Its geodesic flow is generated by
a fiberwise algebraic vector field. Indeed, locally, this vector field has
the form:
$V: (x, p)= (x_1, \ldots, x_n, p_1, \ldots, p_n)
 \in U \times {\bf R}^n \to (p_1, \ldots, p_n, \Sigma_{i j} \Gamma_{i
j}^1(x) p_ip_j,
\ldots, \Sigma_{i j} \Gamma_{i j}^n(x) p_i p_j)$.
(The
$\Gamma_{ij}^k$ are the christoffel symbols).


Now, we introduce a generalization
of geodesic flows as plane fields on
Grassmann bundles.
Let $ \pi: Gr^d(M) \to M$ be the Grassman bundle of
$d$-planes tangent to
$M$. The connection determines a splitting
$T Gr^d(M) = V \bigoplus H$, where $V$ is the vertical and $H$ is
the horizontal  space (given by
$\nabla$). For $p \in Gr_x^d(M)$, $d_p\pi$ maps isomorphically
$H_p$ onto $T_xM$. Let $\tau^d(p)$ be the $d$-plane
contained in $H_p$ which is mapped by $d_p\pi$ to $p $
(as a subspace of $T_xM$). Thus $\tau^d$ is a
$d$-plane field on $Gr^d(M)$, called
the {\it tautological geodesic} plane field on $Gr^d(M)$.
(We think that this construction must be known, although we
haven't found any reference where it is
explicitly mentioned, see \cite{Zeghib.tautologic} for more details
and a systematic study).

The tautological character of $\tau^d$ is clear. The geodesic
adjective is justified by  the fact that, the projection of a leaf
of $\tau^d$ is a (totally)
geodesic submanifold of dimension $d$ in $M$. Conversely, if
$S$ is a $d$-dimensional geodesic
submanifold of $M$, then its Gauss lift $x \in S \to T_xS \in Gr^d(M)$
is a leaf of $\tau^d$.

The  fiberwise algebraic discussion
 on
vector   bundles,  extends
in a straightforward way, to projective   bundles (i.e.
fiber bundles whose fibers are projective spaces...). In
particular, here, as in the case of the geodesic flow,
the tautological plane fields $\tau^d$ are
fiberwise algebraic.

In fact, for the following application, we will immediately come back to
a vector bundle situation.

\paragraph{The pseudo-group of local isometries.}
A (local) isometry or a (local)  affine diffeomorphism is
 a local diffeomorphism
of $M$, which preserves $\nabla$. Equivalently, an affine diffeomophism
is a diffeomorphism which
sends (parameterized) geodesics to
(parameterized) geodesics. One may also define { \it  affine mappings}
as, not necessarily diffeomorphic mappings, sending geodesics to geodesics.

One may naturally construct  a product connection $\nabla \bigoplus
\nabla$
on the product
$M \times M$. If $\nabla$ is the Levi-Civita connection of
a pseudo-Riemannian metric $g$, then $\nabla \bigoplus \nabla$
is the Levi-Civita connection of the product
metric $g \bigoplus g$ (which is the same as the Levi Civita
of the product $g \bigoplus -g$). A curve  $t \to
 (c(t), d(t))$ is (a parameterized) geodesic iff
both of its projections $t \to c(t)$, and $t \to d(t)$
are geodesic in $M$.

Let $f:  U \to V$ be  a  smooth map.  Its graph
$Graph(f)$
is a $n$-submanifold of $M \times M$.  One easily sees,
from the characterization of geodesics in $M \times M$, that
$f$ is an affine mapping, iff   $Graph(f)$  is a (totally)
geodesic submanifold in $M \times M$ (the proof works as
in the case of ${\bf R}^n$).

In particular, local affine mappings  give rise to leaves of
 the tautological geodesic plane field
$\tau^n$ on $Gr^n (M \times M)$.

Let $Gr^*(M \times M)$
consist of $n$-planes which are graphs, that is
$p \in Gr^*_{(x, y)} (M \times M)$, iff, $p $
is a graph of a linear map $T_x M \to T_y M$
(or equivalently, $p$ projects injectively on $T_xM$).

Then, a leaf of $\tau^n$
trough an element $p \in Gr^*(M \times M)$, determines
a local affine mapping.

Observe that  $Gr^*(M \times M)$
is a vector   bundle on $M \times M$, the fiber
over $(x, y)$ being $Hom(T_xM \to T_yM)$.

To get  local affine diffeomorphisms, one
considers  $Gr^{**}(M \times M)$, the set of $n$-planes
transverse to each of the factors $M \times \{.\}$
and $\{.\} \times M$, that is,
$p \in Gr^{**}_{(x, y)} (M \times M)$, iff
$p$ is the graph of an isomorphism $T_xM \to T_yM$.

We have the following interpretation: $(x, y)$ belongs to the projection of
the integrability
domain of $\tau^n$ on $Gr^{**}(M \times M)$, iff, there is a local affine
diffeomorphism sending $x$ to $y$, that is $x$ and $y$ have the same
orbit under
the pseudo-group of local affine diffeomorphisms.

It is easy to see $Gr^{**}(M \times M)$ as the complementary in $Gr^*(M
\times M)$
of a fiberwise algebraic set, and hence in particular, it is an (open) fiberwise
constructible set.


\begin{corollary} \label{Gromov}
(Gromov \cite{Gro}, see also \cite{Ben-Gro} and \cite{Fer}) Let $M$
be an affine manifold.
Suppose
that its pseudo-group of local  affine diffeomorphisms
admits  a dense orbit, then,  it has an  open dense orbit (that is
there is an open dense homogeneous set in $M$).

\end{corollary}

\begin{proof} Apply Theorem \ref{constraints} to $P= \tau^n$ on
$Gr^*(M \times M)$, with
a constraint set ${\cal S}= Gr^{**}(M \times M)$.

Let $x_0 \in M$, be  a point with a dense orbit ${\cal O}_0$ under the
 affine pseudo-group.  The projection
of ${\cal S} \cap {\cal D}$ contains ${\cal O}_0 \times
{\cal O}_0 $. From Theorem \ref{constraints}, the
projection of
${\cal S} \cap {\cal D}$ contains an open dense set
$U$ in $M \times M$.
Let $(x, y) \in U$, then the orbit ${\cal O}_x$ of  $x$ under
the  affine
pseudo-group, contains the open (non-empty) set  $(\{x \} \times M) \cap U$
of $\{x\} \times M$.
 Since the orbit ${\cal O}_0$ is dense, we have ${\cal O}_0 \cap {\cal O}_x
\neq \emptyset$, and hence ${\cal O}_0 = {\cal O}_x$.
But obviously, an orbit with non-empty interior is open, therefore
${\cal O}_0$ is open and dense.


\end{proof}

\subsection{Some comments}
\label{counter.example}

\paragraph{Example.}
Consider on ${\bf R}^n$  a connection $\nabla_ = \nabla^0 + T$, where
$\nabla^0$ is the usual flat connection (that is
$\nabla^0_X Y = D_X Y$), and $T = (T_{ij}^k)$ is a symmetric tensor
$T{\bf R}^n \times T{\bf R}^n \to T{\bf R}^n$. Suppose that
$T$ is flat at $0$, that is, all the partial derivatives of all orders,
of the functions
$T_{ij}^k$,  vanish at $0$.

Consider the tautological geodesic  plane field
$\tau^n$ of ${\bf R}^n \times {\bf R}^n$.
 It is easy to see, that
 $ {\cal D}^\infty_{(0,0)} = Gr_{(0,0)}^n( {\bf R}^n \times {\bf R}^n)$. If
the integrability domain ${\cal D}$ contains ${\cal D}^\infty $,
or more precisely, ${\cal D}^\infty_{(0,0)}$, then, in particular, every
linear map $A: {\bf R}^n \to {\bf R}^n$, will be the derivative
of a local affine (for $({\bf R}^n, \nabla)$) map $F_A$  fixing $0$.
It is easy to see that this implies that $T$ is very special.
Indeed, the existence of non-diffeomorphic affine maps,
leads to vanishing relations of the curvature, not
only at $0$, but also near it.


\paragraph{Other constraints.} In the proof of the above corollary,
 one may add
further constraints of algebraic nature. For example, if $M$ is endowed with
a pseudo-Riemannian metric $g$, then one considers $n-$planes of
$M \times M$, which are isotropic with respect to
the pseudo-Riemannian metric $g \bigoplus -g$ on $M \times M$. The obtained
solutions correspond then to local isometries of $(M, g)$. Observe that
the constraint set here is tangent to $\tau^n$, and it is
in fact fiberwise algebraic (not only
fiberwise constructible) in $Gr^*(M \times M)$, since an isotropic plane
which belongs to  $Gr^*(M \times M)$, must belong to $Gr^{**}(M \times M)$.

Similarly, one may treat  the isometry
pseudo-group of  a unimodular affine structure, and in general,
 any algebraic enrichment  of
the affine structure.

\paragraph{The full Gromov's Theorem.} It is the above corollary of
Gromov's theory that
was used in the celebrated work \cite{{Ben-Fou-Lab}} (and also in
\cite{B-L}).

The full Gromov's theorem,
 that is,  for non-necessarily
topologicaly transitive isometry pseudo-groups,
and for general rigid geometric structures
was
utilized in \cite{D'A}, in the analytic case. As
we have said above,
in our approach,
there are no integrability or structure difficulties
in the analytic case. In fact, \cite{D-G}
contains a direct approach in the analytic case.

It is generally admitted  that there are  no serious difficulties to pass from
affine structures to general rigid (algebraic) structures (see for example
a comment in \cite{Ben-Gro}).

Observe that here, just the idea  of affine structures enriched
with algebraic constraints
allows us to generalize Corollary \ref{Gromov} to a large class
of rigid structures (for example that utilized in the proof
of the main result of \cite{D'A}).

Now, for affine structures with non-necessarily
topologicaly transitive isometry pseudo-group,  the idea of the proof of
Gromov's Theorem, is
 to find a submanifold
in $M \times M$, which, ``essentially'', contains as an open subset,
 the projection
of the infinitesimal integrability domain of $\tau^n$.

\paragraph{Compactification. Singular isometries.} We hope that
 our approach here, provides with  elements leading
to
analyze the non-completeness of the locally
homogeneous open dense set
$U$ in $M$. Indeed, ${\cal D}$ is naturally
compactified by ${\cal D}^\infty$, and
there are sometimes strong evidences (as in
the Anosov case of \cite{Ben-Fou-Lab}) that the
set ${\cal D}^\infty - {\cal D}$ must be empty.

Moreover, $Gr^{**}(M \times M)$ is naturally compactified
by $Gr^n (M \times M)$. The (new) leaves of $\tau^n$ in this latter
space, may be interpreted as singular affine mappings,
and from another point of view, as ``stable laminations''
of (regular) affine mappings.

In fact, compactifications may be defined in the general set-up of control
theory
 of  \S \ref{Control}. Indeed, as ${\bf R}^n$ is projectively compactified by
${\bf R}P^n$ (and not ${\bf R}P^{n-1}$), any vector bundle $N \to B$
with fiber type ${\bf R}^n$ can be (fiberwise) compactified
by a  ``projective'' bundle $\bar{N} \to B$, with fiber type
${\bf R}P^n$. Fiberwise objects on $N$ extend to $\bar{N}$, and it seems
interesting
to interpret them there.

\paragraph{Fiberwise algebraic closure.} let's try
to compare the parts of responsibility behind Corollary \ref{Gromov},
of the integrabilty Theorem \ref{hermann}, and the structure  Theorem
\ref{structure}.

For this, let's consider the following
situation. Take   $G$,  a group
of (global) affine diffeomorphisms of $M$. We have a proper
embedding $(g, x) \in G \times M \to Graph(D_xg) \in Gr^{**}(M
\times M)$. Denote its image by $L$.

The projection of $L$ in $M \times M$
is the union of the graphs of all the elements of $G$, and
$L$ itself is nothing but the union of the Gauss lifts of these graphs. For
this,
let's call $L$ the graph of $G$.

For example, if $G$ is discrete and infinite, then the projection of
$L$ is a countable union of graphs. Therefore, from the structure Theorem,
$L$ is far away from being a fiberwise algebraic set (although
it is closed).

It is thus natural to take the {\it fiberwise algebraic closure}
$\overline{L}^{fib, alg}$ of $L$. The structure Theorem ensures  that
$\overline{L}^{fib, alg}$ has a nice projection.

However, one needs to interpret elements of $\overline{L}^{fib, alg}$, in
other words,
one asks, what properties of elements of $L$ pass to
its fiberwise algebraic closure?

It is the integrability Theorem which answers this question by
stating that, away from a nowhere dense set, the new  elements
of $\overline{L}^{fib, alg}$ are local isometries.

In other words, the integrability Theorem states, essentially, that
in contrast with $L$,  the graph of the local isometry pseudo-group
is a fiberwise algebraic set.  The structure Theorem says that
one has won  a lot from the statement of the integrability Theorem.


\begin{remark} {\em
Similarly to the above embedding,
there is a classical way of breaking
dynamics of $G$, by letting it act
on the frame bundle $P \to M$. To keep
everything elementary, compactify $P$ by seeing
it as an open set in $N$, the vector bundle
with fibers, $N_x = Hom({\bf R}^n \to T_x M)$
($n = dim M$). It is endowed with a principal
$GL(n, {\bf R})$-action.

Suppose that the $G$-action on $M$
is topologically transitive, that is,  it has a dense orbit. Then,
there is an open dense set $U $ of $M$, such that for all
$p \in P$, over $U$, the fiberwise-algebraic closure
of $\overline{G.p}^{fib, alg}$ projects onto $U$. Of course, $GL(n, {\bf R})$
 permutes these fiberwise algebraic closures. The stabilizer in
$GL(n, {\bf R})$ of any closure $\overline{G.p}^{fib, alg}$, may be
identified to
the  $C^\infty$-{\it algebraic hull} of $G$, as introduced in
 \cite{Zim}.

One may define in a  natural way, $C^s$-fiberwise
algebraic sets, for any $ s \geq 0$,  and find $C^s$-algebraic
hulls as defined by Zimmer, for all $s \geq 0$.



}

\end{remark}








\section{Proofs}


\subsection{Sketch of proof of Theorem \protect{\ref{hermann}}}
\label{proof.hermann}

Let $\Delta$ be an involutive  locally finitely
generated distribution on $N$.
At $x \in N$, we denote $\Delta(x)$ the evaluation of $\Delta$
 at $x$.

Let $x_0 \in N$. To construct a leaf of $x_0$,
start with a  vector field $V_0$ of $\Delta$,
non singular at $x_0$, and let $\phi^t$ be its flow.
Suppose that $(\phi^t)^*$ preserves  the evaluation of
$\Delta$ along the orbit $\phi^t(x_0)$, that is
$D_{x_0} \phi^t(\Delta(x_0)) = \Delta(\phi^t(x_0))$.
Take another vector field linearly independent of
$V_0$, and let $\psi^t$ be its flow. Suppose that it satisfies the same
invariance condition, then the surface obtained by saturating the
$\phi^t$-orbit by the flow $\psi^t$ is tangent to $\Delta$.

Retiring the construction,
we would obtain a leaf, if we check the
invariance requirement for all vector fields like  $V_0$.


 Locally,
in some co-ordinates system, we may assume $N = {\bf R} \times {\bf R}^{n-1}$,
and $V_0 = \partial / \partial t$. So, $V_0$ generates a translation flow.

Let
$V_1, \ldots, V_k$ a finite set of generating
vector fields of $\Delta$ near $x_0$.

Since  $\Delta$ is involutive,  $[\partial /\partial t, V_i] = \partial V_i
/\partial t
\in \Delta$ (here we see $V_i$ as vectorial maps on ${\bf R}^n$).
Write: $\partial V_i / \partial t =  \Sigma_{1 \leq j \leq k} a_{i j} V_j$.

So, the problem becomes the following, along the $t$-axis
$(t, 0) \in {\bf R} \times {\bf R}^{n-1}$, we are given
vector fields, $V_1(t), \ldots, V_k(t)$,  and  there are smooth functions
$a_{ij}(t)$, such that
$ \partial V_i /\partial t = \Sigma a_{i j} V_j$. Does this imply that
the space generated by $\{ V_1(t), \ldots, V_k(t) \}$ is independent
of $t$ (i.e. it is parallel along the $t$-axis)?

This is clear in the case $k=1$, that is, if a vector field $V(t)$
satisfies a relation
$\partial V / \partial t = a (t) V(t)$, for $a(t)$
continuous, then $V(t)$ has a parallel direction, and if $V(t)$
 vanishes somewhere
then it vanishes everywhere. Indeed, $e(t) = V / \vert V(t)\vert$
is parallel, where it is defined, that is,  where $V(t) \neq 0$. Now, if
for example $0$ is a  boundary point of the set where $V(t) = 0$, then,
on a semi-open interval, say  $[0, \epsilon[$, we have $V(t) = f(t) e$ ($e
= e(t)$),
and thus, $f(t)$ is a $C^0$ non-trivial solution of the equation
$f^\prime = a(t) f(t)$, with $f(0) = 0$, which is impossible.

Next, in the general case, that is $k >1$, near a generic $t$,
it is possible to write, all the vector fields  $V_1, \ldots, V_k$
as smooth combinations of $r$-linearly independent elements,
say $V_1, \ldots, V_r$.  One then considers the exterior product
$ V(t)= V_1(t) \wedge \ldots \wedge V_r(t)$. It satisfies
(near a generic point)
a relation $\partial V / \partial t = a V$, and is therefore
parallel by the first step.

To finish the proof, it suffices to show that the dimension of the space
generated by the $V_i(t)$ is constant. This dimension
equals the rank of the matrix $X = (x_{ij})_{i \leq k, j \leq n}$, defined by
$V_i = \Sigma_j x_{i j}e_j$, where  $(e_i)_{1 \leq i \leq n}$
is the canonical   basis of ${\bf R}^n$.
We have,  $\partial V_i / \partial t
= \Sigma_j (\partial x_{ij} / \partial t) e_j$. On the other hand,
 $\partial V_i / \partial t
= \Sigma_l  a_{il} V_l = \Sigma_{lj} a_{il} x_{lj}e_j$.
Thus $X$ satisfies the equation (on
$k \times n $-matrices) $X^\prime = A(t) X$,
where $A $ is the $k \times k$ matrix $(a_{ij})$.
 Thus $X(t) = R(t) X(0)$, where
$R(t) $ is the $k \times k$-matrix, resolvent
of the equation on ${\bf R}^k$, $Y^\prime = A(t) Y$
($Y \in {\bf R}^k$).
In particular, the rank of $X(t)$ doesn't depend on $t$.

\subsection{Noetherian properties}

We will deal here with polynomials (with many indeterminates)
on  a ring $R$ which is  $C^0(Y)$, the ring of continuous functions on
a topological space $Y$, or $C^k(Y)$, $0 \leq k \leq \infty$,
 the ring
of $k-$differentiable functions on a subset $Y$ of
 a smooth manifold $B$
(Recall that $f \in C^k(Y)$, means
that $f$ extends locally to an element of
$C^k(B)$).



If $Y^\prime $ is a subset of $Y$, there is a restriction homomorphism
$C^k(Y) \to C^k(Y^\prime)$,
 and by the same way restriction homomorphisms
$C^k(Y) [X_1, \ldots, X_m] \to
C^k(Y^\prime) [X_1, \ldots, X_m]$.

This allows us  to restrict other associated
objects, for example, if $I$ is an ideal
of $C^k(Y) [X_1, \ldots, X_m]$, then its restriction
$I \vert Y^\prime$ is  the ideal of
$C^k(Y^\prime) [X_1, \ldots, X_m]$ generated
by the restriction to
$Y^\prime$ of all the elements
of $I$.

An ideal $I$ of $C^k(Y)$ (or
$C^k(Y)[X_1, \ldots, X_m]$) is locally finitely generated, if
every $x \in Y$ admits a neighborhood $U_x$  such that
$I \vert U_x$ is finitely generated.

\begin{lemma}\label{stationnary}
 Let $I_1 \subset \ldots I_j \subset
\ldots$ be an increasing sequence of
ideals in $C^k(Y)$. Then,  there is an open dense
set $U \subset Y$, over
which, all the ideals are locally finitely
generated (that is $I_i \vert U$ is locally
finitely generated, for all $i$), and the sequence of
ideals is locally
stationary (on $U$).
\end{lemma}

\begin{proof} Let $U_x$ be  a neighborhood of $x$,
 and $I$ an ideal, we say that
$I \vert U_x$ is trivial, if either $I \vert U_x
= 0$, or $I \vert U_x = C^k(U_x)$.

Observe that to ensure the
 existence of $U_x$, such that  $I \vert U_x$ equals
$C^k(U_x)$, it suffices that $I$ contains an element
$f $  such that $f(x) \neq 0$.


Let $U = \{x \in Y / $ there is a neighborhood $U_x$  of $x$, such that,
for all $j$, $I_j \vert U_x$ is trivial $\}$.

By definition, $U$ is open. It is clear, that, over $U$, the sequence of ideals
satisfies the requirements of the lemma. Therefore, it
suffices to show that $U$ is dense.

Firstly, $U$ in nonempty. Indeed, let $j$ be the first integer
such that $I_j  \neq 0$. Then, there is $x \in Y$,
and $f \in I_j$, such that, $f(x) \neq 0$, and hence there
is a neighborhood $U_x$, such that $I_j \vert U_x = C^k(U_x)$, and
thus (by definition of $j$), we have:
$0 = I_1 \vert U_x = \ldots I_{j-1} \vert U_x $, and
$C^k(U_x) = I_i \vert U_x$, for all $i \geq j$, that is
$x \in U$.

To see that $U$ is dense, suppose the contrary, and
consider the open (non-empty set)
$Y -\overline{U}$. Restrict everything to it, and
 conclude, as we have just proved,  that its
corresponding $U$, is
non-empty.
Therefore, there is $x$  in   $Y- \overline{U}$, having
a neighborhood $U_x$ (relative to
$Y- \overline{U}$), such that all
the restrictions
$I_j \vert U_x$ are trivial.
But, since $Y - \overline{U}$ is open in
$Y$, $U_x$ is a neighborhood of $x$  in $Y$, and
therefore, by definition, $x \in U$, which contradicts
our hypothesis.

\end{proof}

\begin{nother} \label{nother}

Let $A = C^k(Y) [X_1, \ldots, X_m]$, and
${\goth a} $ an-$A$-submodule of $A^l$ ($l$ is an integer).
Then, there is an open dense
set $U \subset Y$,
 over which,
${\goth a}$ is locally finitely
generated.

\end{nother}

\begin{proof} Firstly, as in the classical case, it suffices to consider
the case $l= 1$, that is ${\goth a}$ is an ideal of $A$.
The
proof (in this case), then follows, as
for Hilbert's basis  Theorem, that is, if a ring $R$ is Noetherian, then
$R[X_1, \ldots, X_m]$ is also Noetherian.

The (classical) proof of this theorem
is done by induction
 on $N$ (see for example \cite{Lang}).
Let's recall how works the reduction from $R[X_1]$ to $R$.
One associates to
the ideal ${\goth a}$ of $R[X_1]$,  an increasing  sequence  $ I_i$
of ideals of $R$, where $I_i$ is the set of elements
appearing as a leading coefficient of an element of ${\goth  a}$ of degree
$\leq i$. One
then arranges a  finitely  generating set  for ${\goth a}$, if
one knows that the sequence is
stationary, and  has at one's disposal finite generating sets for each
$I_i$ (the number of $i$'s in account is finite).

In our case, from Lemma \ref{stationnary},
the sequence of ideals $I_i$, satisfies the finiteness requirements, after
restricting to an open dense set $Y^\prime$. Therefore
${\goth  a} \vert Y^\prime$ is finitely generated.


\end{proof}



\medskip
\noindent
CNRS, UMPA, \'Ecole Normale Sup\'erieure de Lyon \\
46, all\'ee d'Italie,
 69364 Lyon cedex 07,  FRANCE \\
Zeghib@umpa.ens-lyon.fr \\
http://umpa.ens-lyon.fr/\~{}zeghib/

\end{document}